\documentclass[12pt]{amsart}
\usepackage{amssymb,framed,mathtools,pgfplots,verbatim}
\definecolor{hot}{RGB}{65,105,225} 
\usepackage[colorlinks=true, allcolors=hot]{hyperref}
\usepackage[margin = 1in]{geometry}

\pgfplotsset{compat=1.18}

\usepackage{listings}
\usepackage{xcolor}
\definecolor{maccolor}{rgb}{0.3,0.3,0.8}
\lstdefinelanguage{Macaulay2}
{
basicstyle={\ttfamily},
keywordstyle={\color{maccolor!80!black}},
commentstyle={\color{gray}},
stringstyle={\color{red!40!black}},
rulecolor=\color{maccolor},
basewidth={1.2ex}, 
sensitive=false,
morecomment=[l]{--},
morecomment=[s]{-*}{*-},
morestring=[b]",
escapechar={`},
escapebegin={\rmfamily},
morekeywords={about,abs,AbstractToricVarieties,accumulate,Acknowledgement,acos,acosh,acot,addCancelTask,addDependencyTask,addEndFunction,addHook,AdditionalPaths,addStartFunction,addStartTask,Adjacent,adjoint,AdjointIdeal,AffineVariety,AfterEval,AfterNoPrint,AfterPrint,agm,AInfinity,alarm,AlgebraicSplines,Algorithm,Alignment,all,AllCodimensions,allowableThreads,ambient,analyticSpread,Analyzer,AnalyzeSheafOnP1,ancestor,ancestors,ANCHOR,and,andP,AngleBarList,ann,annihilator,antipode,any,append,applicationDirectory,applicationDirectorySuffix,apply,applyKeys,applyPairs,applyTable,applyValues,apropos,argument,Array,arXiv,Ascending,ascii,asin,asinh,ass,assert,associatedGradedRing,associatedPrimes,AssociativeAlgebras,AssociativeExpression,atan,atan2,atEndOfFile,Authors,autoload,AuxiliaryFiles,backtrace,Bag,Bareiss,baseFilename,BaseFunction,baseName,baseRing,baseRings,BaseRow,BasicList,basis,BasisElementLimit,Bayer,BeforePrint,beginDocumentation,BeginningMacaulay2,Benchmark,benchmark,Bertini,BesselJ,BesselY,betti,BettiCharacters,BettiTally,between,BGG,BIBasis,Binary,BinaryOperation,Binomial,binomial,BinomialEdgeIdeals,Binomials,BKZ,BlockMatrix,BLOCKQUOTE,BODY,Body,BoijSoederberg,BOLD,Book3264Examples,Boolean,BooleanGB,borel,Boxes,BR,break,Browse,Bruns,cache,CacheExampleOutput,CacheFunction,CacheTable,cacheValue,CallLimit,cancelTask,capture,catch,Caveat,CC,CDATA,ceiling,Center,centerString,Certification,ChainComplex,chainComplex,ChainComplexExtras,ChainComplexMap,ChainComplexOperations,ChangeMatrix,char,CharacteristicClasses,characters,charAnalyzer,check,CheckDocumentation,chi,Chordal,class,Classic,clean,clearAll,clearEcho,clearOutput,close,closeIn,closeOut,ClosestFit,CODE,code,codim,CodimensionLimit,coefficient,CoefficientRing,coefficientRing,coefficients,Cofactor,CohenEngine,CohenTopLevel,CoherentSheaf,CohomCalg,cohomology,coimage,CoincidentRootLoci,coker,cokernel,collectGarbage,columnAdd,columnate,columnMult,columnPermute,columnRankProfile,columnSwap,combine,Command,commandInterpreter,commandLine,COMMENT,commonest,commonRing,comodule,CompactMatrix,compactMatrixForm,CompiledFunction,CompiledFunctionBody,CompiledFunctionClosure,Complement,complement,complete,CompleteIntersection,CompleteIntersectionResolutions,Complexes,ComplexField,components,compose,compositions,compress,concatenate,conductor,ConductorElement,cone,Configuration,ConformalBlocks,conjugate,connectionCount,Consequences,Constant,Constants,constParser,content,continue,contract,Contributors,ConvexInterface,conwayPolynomial,ConwayPolynomials,copy,copyDirectory,copyFile,copyright,Core,CorrespondenceScrolls,cos,cosh,cot,CotangentSchubert,cotangentSheaf,coth,cover,coverMap,cpuTime,createTask,Cremona,csc,csch,current,currentColumnNumber,currentDirectory,currentFileDirectory,currentFileName,currentLayout,currentLineNumber,currentPackage,currentString,currentTime,Cyclotomic,Database,Date,DD,dd,deadParser,debug,debugError,DebuggingMode,debuggingMode,debugLevel,DecomposableSparseSystems,Decompose,decompose,deepSplice,Default,default,defaultPrecision,Degree,degree,degreeLength,DegreeLift,DegreeLimit,DegreeMap,DegreeOrder,DegreeRank,Degrees,degrees,degreesMonoid,degreesRing,delete,demark,denominator,Dense,Density,Depth,depth,Descending,Descent,Describe,describe,Description,det,determinant,DeterminantalRepresentations,DGAlgebras,diagonalMatrix,diameter,Dictionary,dictionary,dictionaryPath,diff,DiffAlg,difference,dim,directSum,disassemble,discriminant,dismiss,Dispatch,distinguished,DIV,Divide,divideByVariable,DivideConquer,DividedPowers,Divisor,DL,Dmodules,do,doc,docExample,docTemplate,document,DocumentTag,Down,drop,DT,dual,eagonNorthcott,EagonResolution,echoOff,echoOn,EdgeIdeals,edit,EigenSolver,eigenvalues,eigenvectors,eint,EisenbudHunekeVasconcelos,elapsedTime,elapsedTiming,elements,Eliminate,eliminate,Elimination,EliminationMatrices,EllipticCurves,EllipticIntegrals,else,EM,Email,End,end,endl,endPackage,Engine,engineDebugLevel,EngineRing,EngineTests,entries,EnumerationCurves,environment,Equation,EquivariantGB,erase,erf,erfc,error,errorDepth,euler,EulerConstant,eulers,even,EXAMPLE,ExampleFiles,ExampleItem,examples,ExampleSystems,Exclude,exec,exit,exp,expectedReesIdeal,expm1,exponents,export,exportFrom,exportMutable,Expression,expression,Ext,extend,ExteriorIdeals,ExteriorModules,exteriorPower,Factor,factor,false,Fano,FastMinors,FastNonminimal,FGLM,File,fileDictionaries,fileExecutable,fileExists,fileExitHooks,fileLength,fileMode,FileName,FilePosition,fileReadable,fileTime,fileWritable,fillMatrix,findFiles,findHeft,FindOne,findProgram,findSynonyms,FiniteFittingIdeals,First,first,firstkey,FirstPackage,fittingIdeal,flagLookup,FlatMonoid,flatten,flattenRing,Flexible,flip,floor,flush,fold,FollowLinks,for,forceGB,fork,FormalGroupLaws,Format,format,formation,FourierMotzkin,FourTiTwo,fpLLL,frac,fraction,FractionField,frames,FrobeniusThresholds,from,fromDividedPowers,fromDual,Function,FunctionApplication,FunctionBody,functionBody,FunctionClosure,FunctionFieldDesingularization,fusePairs,futureParser,GaloisField,gb,GBDegrees,gbRemove,gbSnapshot,gbTrace,gcd,gcdCoefficients,gcdLLL,GCstats,genera,GeneralOrderedMonoid,GenerateAssertions,generateAssertions,generator,generators,Generic,GenericInitialIdeal,genericMatrix,genericSkewMatrix,genericSymmetricMatrix,gens,genus,get,getc,getChangeMatrix,getenv,getGlobalSymbol,getNetFile,getNonUnit,getPrimeWithRootOfUnity,getSymbol,getWWW,GF,gfanInterface,Givens,GKMVarieties,GLex,Global,global,globalAssign,globalAssignFunction,GlobalAssignHook,globalAssignment,globalAssignmentHooks,GlobalDictionary,GlobalHookStore,globalReleaseFunction,GlobalReleaseHook,Gorenstein,GradedLieAlgebras,GradedModule,gradedModule,GradedModuleMap,gradedModuleMap,gramm,GraphicalModels,GraphicalModelsMLE,Graphics,graphIdeal,graphRing,Graphs,Grassmannian,GRevLex,GroebnerBasis,groebnerBasis,GroebnerBasisOptions,GroebnerStrata,GroebnerWalk,groupID,GroupLex,GroupRevLex,GTZ,Hadamard,handleInterrupts,HardDegreeLimit,hash,HashTable,hashTable,HEAD,HEADER1,HEADER2,HEADER3,HEADER4,HEADER5,HEADER6,HeaderType,Heading,Headline,Heft,heft,Height,height,help,Hermite,hermite,Hermitian,HH,hh,HigherCIOperators,HighestWeights,Hilbert,hilbertFunction,hilbertPolynomial,hilbertSeries,HodgeIntegrals,hold,Holder,Hom,homeDirectory,HomePage,Homogeneous,Homogeneous2,homogenize,homology,homomorphism,HomotopyLieAlgebra,hooks,horizontalJoin,HorizontalSpace,HR,HREF,HTML,html,httpHeaders,Hybrid,HyperplaneArrangements,Hypertext,hypertext,HypertextContainer,HypertextParagraph,icFracP,icFractions,icMap,icPIdeal,id,Ideal,ideal,idealizer,identity,if,IgnoreExampleErrors,ii,image,imaginaryPart,IMG,ImmutableType,importFrom,in,incomparable,Increment,independentSets,indeterminate,IndeterminateNumber,Index,index,indexComponents,IndexedVariable,IndexedVariableTable,indices,inducedMap,inducesWellDefinedMap,InexactField,InexactFieldFamily,InexactNumber,InfiniteNumber,infinity,info,InfoDirSection,infoHelp,Inhomogeneous,input,Inputs,insert,installAssignmentMethod,installedPackages,installHilbertFunction,installMethod,installMinprimes,installPackage,InstallPrefix,instance,instances,IntegralClosure,integralClosure,integrate,IntermediateMarkUpType,interpreterDepth,intersect,intersectInP,Intersection,intersection,interval,InvariantRing,inverse,InverseMethod,inversePermutation,Inverses,inverseSystem,InverseSystems,Invertible,InvolutiveBases,irreducibleCharacteristicSeries,irreducibleDecomposition,isAffineRing,isANumber,isBorel,isCanceled,isCommutative,isConstant,isDirectory,isDirectSum,isEmpty,isField,isFinite,isFinitePrimeField,isFreeModule,isGlobalSymbol,isHomogeneous,isIdeal,isInfinite,isInjective,isInputFile,isIsomorphism,isLinearType,isListener,isLLL,isMember,isModule,isMonomialIdeal,isNormal,isOpen,isOutputFile,isPolynomialRing,isPrimary,isPrime,isPrimitive,isPseudoprime,isQuotientModule,isQuotientOf,isQuotientRing,isReady,isReal,isReduction,isRegularFile,isRing,isSkewCommutative,isSorted,isSquareFree,isStandardGradedPolynomialRing,isSubmodule,isSubquotient,isSubset,isSupportedInZeroLocus,isSurjective,isTable,isUnit,isWellDefined,isWeylAlgebra,ITALIC,Iterate,Jacobian,jacobian,jacobianDual,Jets,Join,join,Jupyter,K3Carpets,K3Surfaces,Keep,KeepFiles,KeepZeroes,ker,kernel,kernelLLL,kernelOfLocalization,Key,keys,Keyword,Keywords,kill,koszul,Kronecker,KustinMiller,LABEL,last,lastMatch,LATER,LatticePolytopes,Layout,lcm,leadCoefficient,leadComponent,leadMonomial,leadTerm,Left,left,length,LengthLimit,letterParser,Lex,LexIdeals,LI,Licenses,LieTypes,lift,liftable,Limit,limitFiles,limitProcesses,Linear,LinearAlgebra,LinearTruncations,lineNumber,lines,LINK,linkFile,List,list,listForm,listLocalSymbols,listSymbols,listUserSymbols,LITERAL,LLL,LLLBases,lngamma,load,loadDepth,LoadDocumentation,loadedFiles,loadedPackages,loadPackage,Local,local,localDictionaries,LocalDictionary,localize,LocalRings,locate,log,log1p,LongPolynomial,lookup,lookupCount,LowerBound,LUdecomposition,M0nbar,M2CODE,Macaulay2Doc,makeDirectory,MakeDocumentation,makeDocumentTag,MakeHTML,MakeInfo,MakeLinks,makePackageIndex,MakePDF,makeS2,Manipulator,map,MapExpression,MapleInterface,markedGB,Markov,MarkUpType,match,mathML,Matrix,matrix,MatrixExpression,Matroids,max,maxAllowableThreads,maxExponent,MaximalRank,maxPosition,MaxReductionCount,MCMApproximations,member,memoize,memoizeClear,memoizeValues,MENU,merge,mergePairs,META,method,MethodFunction,MethodFunctionBinary,MethodFunctionSingle,MethodFunctionWithOptions,methodOptions,methods,midpoint,min,minExponent,mingens,mingle,minimalBetti,MinimalGenerators,MinimalMatrix,minimalPresentation,minimalPresentationMap,minimalPresentationMapInv,MinimalPrimes,minimalPrimes,minimalReduction,Minimize,minimizeFilename,MinimumVersion,minors,minPosition,minPres,minprimes,Minus,minus,Miura,MixedMultiplicity,mkdir,mod,Module,module,ModuleDeformations,modulo,MonodromySolver,Monoid,monoid,MonoidElement,Monomial,MonomialAlgebras,monomialCurveIdeal,MonomialIdeal,monomialIdeal,MonomialIntegerPrograms,MonomialOrbits,MonomialOrder,Monomials,monomials,MonomialSize,monomialSubideal,moveFile,multidegree,multidoc,multigraded,MultigradedBettiTally,MultiGradedRationalMap,multiplicity,MultiplicitySequence,MultiplierIdeals,MultiplierIdealsDim2,MultiprojectiveVarieties,mutable,MutableHashTable,mutableIdentity,MutableList,MutableMatrix,mutableMatrix,NAGtypes,Name,nanosleep,Nauty,NautyGraphs,NCAlgebra,NCLex,needs,needsPackage,Net,net,NetFile,netList,new,newClass,newCoordinateSystem,NewFromMethod,newline,NewMethod,newNetFile,NewOfFromMethod,NewOfMethod,newPackage,newRing,nextkey,nextPrime,nil,NNParser,NoetherianOperators,NoetherNormalization,NonminimalComplexes,nonspaceAnalyzer,NoPrint,norm,normalCone,Normaliz,NormalToricVarieties,not,Nothing,notify,notImplemented,NTL,null,nullaryMethods,nullhomotopy,nullParser,nullSpace,Number,number,NumberedVerticalList,numcols,numColumns,numerator,numeric,NumericalAlgebraicGeometry,NumericalCertification,NumericalImplicitization,NumericalLinearAlgebra,NumericalSchubertCalculus,numericInterval,NumericSolutions,numgens,numRows,numrows,odd,oeis,of,ofClass,OL,OldPolyhedra,OldToricVectorBundles,on,OneExpression,OnlineLookup,OO,oo,ooo,oooo,openDatabase,openDatabaseOut,openFiles,openIn,openInOut,openListener,OpenMath,openOut,openOutAppend,operatorAttributes,Option,OptionalComponentsPresent,optionalSignParser,Options,options,OptionTable,optP,or,Order,order,OrderedMonoid,orP,OutputDictionary,Outputs,override,pack,Package,package,PackageCitations,PackageDictionary,PackageExports,PackageImports,PackageTemplate,packageTemplate,pad,pager,PairLimit,pairs,PairsRemaining,PARA,Parametrization,parent,Parenthesize,Parser,Parsing,part,Partition,partition,partitions,parts,path,pdim,peek,PencilsOfQuadrics,Permanents,permanents,permutations,pfaffians,PHCpack,PhylogeneticTrees,pi,PieriMaps,pivots,PlaneCurveSingularities,plus,poincare,poincareN,Points,polarize,poly,Polyhedra,Polymake,PolynomialRing,Posets,Position,position,positions,PositivityToricBundles,POSIX,Postfix,Power,power,powermod,PRE,Precision,precision,Prefix,prefixDirectory,prefixPath,preimage,prepend,presentation,pretty,primaryComponent,PrimaryDecomposition,primaryDecomposition,PrimaryTag,PrimitiveElement,Print,print,printerr,printingAccuracy,printingLeadLimit,printingPrecision,printingSeparator,printingTimeLimit,printingTrailLimit,printString,printWidth,processID,Product,product,ProductOrder,profile,profileSummary,Program,programPaths,ProgramRun,Proj,Projective,ProjectiveHilbertPolynomial,projectiveHilbertPolynomial,ProjectiveVariety,promote,protect,Prune,prune,PruneComplex,pruningMap,Pseudocode,pseudocode,pseudoRemainder,Pullback,PushForward,pushForward,Python,QQ,QQParser,QRDecomposition,QthPower,Quasidegrees,QuaternaryQuartics,QuillenSuslin,quit,Quotient,quotient,quotientRemainder,QuotientRing,Radical,radical,RadicalCodim1,radicalContainment,RaiseError,random,RandomCanonicalCurves,RandomComplexes,RandomCurves,RandomCurvesOverVerySmallFiniteFields,RandomGenus14Curves,RandomIdeals,randomKRationalPoint,RandomMonomialIdeals,randomMutableMatrix,RandomObjects,RandomPlaneCurves,RandomPoints,RandomSpaceCurves,Range,rank,RationalMaps,RationalPoints,RationalPoints2,ReactionNetworks,read,readDirectory,readlink,readPackage,RealField,RealFP,realPart,realpath,RealQP,RealQP1,RealRoots,RealRR,RealXD,recursionDepth,recursionLimit,Reduce,reducedRowEchelonForm,reduceHilbert,reductionNumber,ReesAlgebra,reesAlgebra,reesAlgebraIdeal,reesIdeal,References,ReflexivePolytopesDB,regex,regexQuote,registerFinalizer,regSeqInIdeal,Regularity,regularity,relations,RelativeCanonicalResolution,relativizeFilename,Reload,remainder,RemakeAllDocumentation,remove,removeDirectory,removeFile,removeLowestDimension,reorganize,replace,RerunExamples,res,reshape,ResidualIntersections,ResLengthThree,Resolution,resolution,ResolutionsOfStanleyReisnerRings,restart,Result,resultant,Resultants,return,returnCode,Reverse,reverse,RevLex,Right,right,Ring,ring,RingElement,RingFamily,ringFromFractions,RingMap,rootPath,roots,rootURI,rotate,round,rowAdd,RowExpression,rowMult,rowPermute,rowRankProfile,rowSwap,RR,RRi,rsort,run,RunDirectory,RunExamples,RunExternalM2,runHooks,runLengthEncode,runProgram,same,saturate,Saturation,scan,scanKeys,scanLines,scanPairs,scanValues,schedule,schreyerOrder,Schubert,Schubert2,SchurComplexes,SchurFunctors,SchurRings,SCRIPT,scriptCommandLine,ScriptedFunctor,SCSCP,searchPath,sec,sech,SectionRing,SeeAlso,seeParsing,SegreClasses,select,selectInSubring,selectVariables,SelfInitializingType,SemidefiniteProgramming,Seminormalization,separate,SeparateExec,separateRegexp,Sequence,sequence,Serialization,serialNumber,Set,set,setEcho,setGroupID,setIOExclusive,setIOSynchronized,setIOUnSynchronized,setRandomSeed,setup,setupEmacs,sheaf,SheafExpression,sheafExt,sheafHom,SheafOfRings,shield,ShimoyamaYokoyama,short,show,showClassStructure,showHtml,showStructure,showTex,showUserStructure,SimpleDoc,simpleDocFrob,SimplicialComplexes,SimplicialDecomposability,SimplicialPosets,SimplifyFractions,sin,singularLocus,sinh,size,size2,SizeLimit,SkewCommutative,SlackIdeals,sleep,SLnEquivariantMatrices,SLPexpressions,SMALL,smithNormalForm,solve,someTerms,Sort,sort,sortColumns,SortStrategy,source,SourceCode,SourceRing,SPACE,SpaceCurves,SPAN,span,SparseMonomialVectorExpression,SparseResultants,SparseVectorExpression,Spec,SpechtModule,SpecialFanoFourfolds,specialFiber,specialFiberIdeal,SpectralSequences,splice,splitWWW,sqrt,SRdeformations,stack,stacksProject,Standard,standardForm,standardPairs,StartWithOneMinor,stashValue,StatePolytope,StatGraphs,status,stderr,stdio,step,StopBeforeComputation,stopIfError,StopWithMinimalGenerators,Strategy,String,STRONG,StronglyStableIdeals,STYLE,Style,style,SUB,sub,SubalgebraBases,sublists,submatrix,submatrixByDegrees,Subnodes,subquotient,SubringLimit,Subscript,subscript,SUBSECTION,subsets,substitute,substring,subtable,Sugarless,Sum,sum,SumOfTwists,SumsOfSquares,SUP,super,SuperLinearAlgebra,Superscript,superscript,support,SVD,SVDComplexes,switch,SwitchingFields,sylvesterMatrix,Symbol,symbol,SymbolBody,symbolBody,SymbolicPowers,symlinkDirectory,symlinkFile,symmetricAlgebra,symmetricAlgebraIdeal,symmetricKernel,SymmetricPolynomials,symmetricPower,synonym,SYNOPSIS,syz,Syzygies,SyzygyLimit,SyzygyMatrix,SyzygyRows,syzygyScheme,TABLE,Table,table,take,Tally,tally,tan,TangentCone,tangentCone,tangentSheaf,tanh,target,Task,taskResult,TateOnProducts,TD,temporaryFileName,tensor,tensorAssociativity,TensorComplexes,terminalParser,terms,TEST,Test,testExample,testHunekeQuestion,TestIdeals,TestInput,tests,TEX,tex,TeXmacs,texMath,Text,TH,then,Thing,ThinSincereQuivers,ThreadedGB,threadVariable,Threshold,throw,Time,time,times,timing,TITLE,TO,to,TO2,toAbsolutePath,toCC,toDividedPowers,toDual,toExternalString,toField,TOH,toList,toLower,top,top,topCoefficients,Topcom,topComponents,topLevelMode,Tor,TorAlgebra,Toric,ToricInvariants,ToricTopology,ToricVectorBundles,toRR,toRRi,toSequence,toString,TotalPairs,toUpper,TR,trace,transpose,TriangularSets,Tries,Trim,trim,Triplets,Tropical,true,Truncate,truncate,truncateOutput,Truncations,try,TSpreadIdeals,TT,tutorial,Type,TypicalValue,typicalValues,UL,ultimate,unbag,uncurry,Undo,undocumented,uniform,uninstallAllPackages,uninstallPackage,Unique,unique,Units,Unmixed,unsequence,unstack,Up,UpdateOnly,UpperTriangular,URL,urlEncode,Usage,use,UseCachedExampleOutput,UseHilbertFunction,UserMode,userSymbols,UseSyzygies,utf8,utf8check,validate,value,values,Variable,VariableBaseName,Variables,Variety,variety,vars,Vasconcelos,Vector,vector,VectorExpression,VectorFields,VectorGraphics,Verbose,Verbosity,Verify,VersalDeformations,versalEmbedding,Version,version,VerticalList,VerticalSpace,viewHelp,VirtualResolutions,VirtualTally,VisibleList,Visualize,wait,WebApp,wedgeProduct,weightRange,Weights,WeylAlgebra,WeylGroups,when,whichGm,while,width,wikipedia,Wrap,wrap,WrapperType,XML,xor,youngest,zero,ZeroExpression,zeta,ZZ,ZZParser,
makeGWClass,getDiagonalClass,makeDiagonalForm,getSignature,isAnisotropic,isIsotropic,getAnisotropicPart,getSumDecompositionString,getGlobalA1Degree,getLocalA1Degree,isIsomorphicForm,addGW,getSumDecomposition}
}
\lstalias{Macaulay2output}{Macaulay2}

\newtheorem{theorem}{Theorem}

\newtheorem{proposition}[theorem]{Proposition}

\theoremstyle{definition}
\newtheorem{definition}[theorem]{Definition}
\newtheorem{example}[theorem]{Example}
\newtheorem{remark}[theorem]{Remark}

\providecommand{\A}{\mathbb{A}}

\DeclareMathOperator{\Gr}{Gr}
\DeclareMathOperator{\GW}{GW}

\DeclareMathOperator{\Spec}{Spec}
\DeclareMathOperator{\Sym}{Sym}

\title{$\A^1$-Brouwer degrees in Macaulay2}

\author{Nikita Borisov} 
\address[N.~Borisov]{Department of Mathematics\\
         University of Pennsylvania\\ 
         Philadelphia, PA 19104\\ 
         USA} 
\email{nborisov@sas.upenn.edu}
\urladdr{https://nikita-borisov.github.io/}
\author{Thomas Brazelton} 
\address[T.~Brazelton]{Department of Mathematics\\ 
         Harvard University\\ 
         Cambridge, MA 02138\\
         USA} 
\email{brazelton@math.harvard.edu}
\urladdr{https://tbrazel.github.io/}
\author{Frenly Espino} 
\address[F.~Espino]{Department of Mathematics\\
         University of Pennsylvania\\ 
         Philadelphia, PA 19104\\ 
         USA} 
\email{frenly@sas.upenn.edu }
\urladdr{https://frenlye.github.io/}
\author{Thomas Hagedorn} 
\address[T.~Hagedorn]{Department of Mathematics and Statistics\\ 
         The College of New Jersey\\ 
         Ewing, NJ 08618\\ 
         USA} 
\email{hagedorn@tcnj.edu}
\urladdr{https://hagedorn.pages.tcnj.edu/}
\author{Zhaobo Han} 
\address[Z.~Han]{University of Pennsylvania\\ 
         Philadelphia, PA 19104\\ 
         USA} 
\email{zbtomhan@sas.upenn.edu}
\urladdr{https://www.linkedin.com/in/zhaobo-han-77b1301a2/}
\author{Jordy Lopez Garcia} 
\address[J.~Lopez Garcia]{Department of Mathematics\\ 
         Texas A\&M University\\ 
         College Station, TX 77843\\ 
         USA} 
\email{jordy.lopez@tamu.edu}
\urladdr{https://jordylopez27.github.io/}
\author{Joel Louwsma} 
\address[J.~Louwsma]{Department of Mathematics\\ 
         Niagara University\\ 
         Niagara University, NY 14109\\ 
         USA} 
\email{jlouwsma@niagara.edu}
\urladdr{https://www.joellouwsma.com}
\author{Wern Juin Gabriel Ong} 
\address[W.~J.~G.~Ong]{Bowdoin College\\
         Brunswick, ME 04011\\
         USA}
\email{gong@bowdoin.edu}
\urladdr{https://wgabrielong.github.io/}
\author{Andrew R. Tawfeek} 
\address[A.~R.~Tawfeek]{Department of Mathematics\\ 
         University of Washington\\ 
         Seattle, WA 98195\\
         USA} 
\email{atawfeek@uw.edu}
\urladdr{https://atawfeek.com/}


\begin{document}

\begin{abstract}
We describe the \emph{Macaulay2} package \texttt{A1BrouwerDegrees} for computing local and global $\mathbb{A}^1$-Brouwer degrees and studying symmetric bilinear forms over the complex numbers, the real numbers, the rational numbers, and finite fields of characteristic not equal to~$2$. 
\end{abstract}

\maketitle

\section{Introduction}

In $\A^1$-homotopy theory, the \emph{$\A^1$-Brouwer degree} provides an algebro-geometric analogue of the classical Brouwer degree from differential topology. Morel's $\A^1$-degree homomorphism identifies the zeroth stable stem of the motivic sphere spectrum with the \emph{Grothendieck--Witt ring} of symmetric bilinear forms over a field \cite[Corollary~1.24]{Morel}. Given an endomorphism of affine space with an isolated rational zero, work of Kass and Wickelgren~\cite{KW-EKL} identifies its local $\A^1$-Brouwer degree with the Eisenbud--Khimshiashvili--Levine signature form \cite{EL,Khim}, which was used to compute local Brouwer degrees in real differential topology. Work of Bachmann and Wickelgren~\cite{BW} extends this work, identifying the $\A^1$-Brouwer degree with a quadratic Grothendieck--Serre duality form.

In \emph{$\A^1$-enumerative geometry} \cite{CubicKW,Levine} (see \cite{Braz-expos,Pauli-expos} for an overview), the $\A^1$-Brouwer degree has found a wealth of applications, recently including \cite{AK,KSLW,JPP}. For instance, via McKean's B\'ezout theorem, the $\A^1$-Brouwer degree can be understood as a quadratically enriched analogue of intersection multiplicity, often encoding deeper geometric information than was available over the algebraically closed fields -- with other invariants of the quadratic form over~$k$ capturing field-specific arithmetic data~\cite{McKean}.

Recent work of the second-named author, McKean, and Pauli~\cite{BMP} provides tractable formulas for computing $\A^1$-Brouwer degrees as \emph{B\'ezoutian bilinear forms}. In the \texttt{A1BrouwerDegrees} package, we implement these methods in \emph{Macaulay2}~\cite{M2} over the fields $\mathbb{C}$, $\mathbb{R}$, $\mathbb{Q}$, and~$\mathbb{F}_q$ (for $q$ odd) and provide a suite of tools whose capabilities include:
\begin{enumerate}
    \item computing $\A^1$-Brouwer degrees (both local and global) for endomorphisms of affine space;\footnote{Due to~$\mathbb{R}$ being an inexact field, $\mathbb{A}^{1}$-Brouwer degrees over~$\mathbb{R}$ have to be computed over~$\mathbb{Q}$ and base-changed to~$\mathbb{R}$.}
    \item decomposing symmetric bilinear forms into their isotropic and anisotropic parts; and
    \item extracting invariants of symmetric bilinear forms (rank, signature, discriminant, Hasse--Witt invariants).
\end{enumerate}

\begin{remark} 
The current scope of this package is the complex numbers, the real numbers, the rational numbers, and finite fields of characteristic not equal to~$2$, so all of our algorithms are implicitly taking place in these settings. We hope to expand the scope of these algorithms in future work.
\end{remark}

In Section~\ref{sec:GWring}, we provide a rapid introduction to the theory of symmetric bilinear forms, highlighting the capacity of our package to build forms, check isomorphisms, and decompose forms. In Section~\ref{sec:A1Degs}, we discuss local and global $\A^1$-Brouwer degrees and provide some computational examples, including quadratically enriched intersection multiplicity of real curves, the $\A^1$-Euler characteristic of the Grassmannian $\Gr(2,4)$ (following \cite[Section~8.2]{BMP}), and local computations for $27$~lines on a cubic surface (following \cite{CubicKW,Pauli}).

\subsection{Software availability}
The software documented here is available in versions 1.23 and later of \emph{Macaulay2} as the \href{https://macaulay2.com/doc/Macaulay2/share/doc/Macaulay2/A1BrouwerDegrees/html/index.html}{\texttt{A1BrouwerDegrees} package}.

\section{The Grothendieck--Witt ring}\label{sec:GWring}

For this entire section, we assume $k$ is a field of characteristic not equal to~$2$. We say a bilinear form $\beta \colon V \times V \to k$ is \emph{symmetric} if $\beta(v,w) = \beta(w,v)$ for all $v,w\in V$. We say $\beta$ is \emph{non-degenerate} if $\beta(v,-) \colon V \to k$ is identically zero if and only if $v=0$.

\begin{definition}\label{def:gram-matrix} 
Let $\beta \colon V \times V \to k$ be a symmetric bilinear form, and choose a basis $e_1, \dots, e_n$ for~$V$. We define the \emph{Gram matrix} of~$\beta$ in the basis $\{e_i\}_{i=1}^n$ to be the symmetric matrix with entries $\beta(e_i,e_j)$. 
\end{definition}

\begin{remark}
    Non-degeneracy of~$\beta$ is equivalent to the statement that the determinant of a Gram matrix in any basis is nonzero. A change of basis for~$V$ corresponds to the associated Gram matrices being congruent.
\end{remark}

Given two symmetric bilinear forms $\beta_i \colon V_i \times V_i \to k$ for $i=1,2$, we can define their sum and product:
\begin{equation}\label{eqn:add-multiply-symm-bilin-forms}
\begin{aligned}
    (\beta_1 \oplus \beta_2) \colon (V_1 \oplus V_2 )\times (V_1 \oplus V_2) &\to k \\
    (\beta_1 \otimes \beta_2) \colon (V_1 \otimes V_2) \times( V_1 \otimes V_2) &\to k.
\end{aligned}
\end{equation}
On Gram matrices, these operations are given by direct sum and tensor product, respectively.

\begin{definition}\label{def:GWk} 
The \emph{Grothendieck--Witt ring} $\GW(k)$ is the group completion of the semiring of isomorphism classes of non-degenerate symmetric bilinear forms over~$k$. 
\end{definition}

\begin{example}
\label{example-5}
Any non-degenerate symmetric bilinear form over an algebraically closed field admits a basis in which its Gram matrix is the identity; therefore rank determines an isomorphism $\GW(\mathbb{C}) \cong \mathbb{Z}$. For further computations of Grothendieck--Witt rings, we refer the reader to \cite[Chapter~II]{Lam}.
\end{example}

When the field~$k$ is the complex numbers, the real numbers, the rational numbers, or a finite field (of characteristic not~$2$), we define a type called \texttt{GrothendieckWittClass} that encodes the class $[\beta]\in \GW(k)$ of a symmetric bilinear form~$\beta$. Grothendieck--Witt classes can be constructed from Gram matrices via the \texttt{makeGWClass} method.

\begin{lstlisting}[language=Macaulay2output]
`\underline{\tt i1}` : needsPackage "A1BrouwerDegrees";
`\underline{\tt i2}` : M = matrix(QQ, {{1,3},{3,7}});
`\underline{\tt o2}` : `$\texttt{Matrix}$ ${\mathbb Q}^{2}\,\longleftarrow \,{\mathbb Q}^{2}$`
`\underline{\tt i3}` : beta = makeGWClass M
`\underline{\tt o3}` = `$\left(\!\begin{array}{cc}
1&3\\
3&7
\end{array}\!\right)$`
`\underline{\tt o3}` : `$\texttt{GrothendieckWittClass}$`
\end{lstlisting}

Given a Grothendieck--Witt class \texttt{beta}, its underlying field can be obtained by running \texttt{getBaseField beta}, and the underlying matrix can be obtained by running either \texttt{beta.matrix} or \texttt{getMatrix beta}. Objects of type \texttt{GrothendieckWittClass} can be added and multiplied via the \texttt{addGW} and \texttt{multiplyGW} methods, respectively.

\begin{example}\label{ex:labelname} 
For any unit $a\in k^\times$, there is a symmetric bilinear form of rank one
\begin{align*}
    \langle a \rangle \colon k \times k &\to k \\
    (x,y) &\mapsto axy.
\end{align*}
Via the change of basis $(x,y)\mapsto (bx,by)$ for any unit $b\in k^\times$, we observe that $\langle a \rangle = \langle ab^2 \rangle$. Hence the representative for a rank one form is determined only by its square class.
\end{example}

The following classical result (see \cite[Corollary I.2.4]{Lam}) implies that the forms $\langle a \rangle$ generate $\GW(k)$.

\begin{theorem}\label{thm:sylvester} 
Every class in $\GW(k)$ is represented by a diagonal Gram matrix.
\end{theorem}

A diagonal representative of a class in $\GW(k)$ can be found using the \texttt{getDiagonalClass} method. 

\begin{lstlisting}[language=Macaulay2output]
`\underline{\tt i4}` : getDiagonalClass beta
`\underline{\tt o4}` = `$\left(\!\begin{array}{cc}
1&0\\
0&-2
\end{array}\!\right)$`
`\underline{\tt o4}` : `$\texttt{GrothendieckWittClass}$`
\end{lstlisting}

We also provide methods for constructing various forms. We can construct a class corresponding to a list of diagonal entries via the \texttt{makeDiagonalForm} method. We denote by $\langle a_1,\ldots,a_n\rangle$ the direct sum of the rank one forms $\langle a_i\rangle$ for $1\le i \le n$.

\begin{lstlisting}[language=Macaulay2output]
`\underline{\tt i5}` : makeDiagonalForm(GF(13), (2,6))
`\underline{\tt o5}` = `$\left(\!\begin{array}{cc}
2&0\\
0&6
\end{array}\!\right)$`
`\underline{\tt o5}` : `$\texttt{GrothendieckWittClass}$`
\end{lstlisting}

Hyperbolic forms are crucial objects of study due to their local-to-global behavior (see Theorem~\ref{thm:lgp-isotropy}), and they can be produced via the \texttt{makeHyperbolicForm} method. Similarly, Pfister forms, which are important objects of study in the world of quadratic forms \cite[Chapter~X]{Lam}, can be produced via the \texttt{makePfisterForm} method.

\subsection{Verifying isomorphisms of forms}\label{subsec: isomorphisms of forms}
Given two non-degenerate symmetric bilinear forms, a natural question is whether they represent the same element of $\GW(k)$. An easy invariant to check is whether they are defined on vector spaces of the same dimension, i.e., whether the \emph{rank} of the forms (the rank of their Gram matrices) agrees. As mentioned in Example~\ref{example-5}, since $\mathbb{C}$ is algebraically closed and every number is a square, rank completely classifies symmetric bilinear forms over the complex numbers. 

Since there are two square classes over the real numbers, namely $+1$ and~$-1$, we can find a Gram matrix representative of any form which is diagonal, with only $\pm1$ appearing along the diagonal. The trace of such a Gram matrix is an invariant of the form, called the \emph{signature}. Rank and signature jointly classify symmetric bilinear forms over the real numbers.

\begin{lstlisting}[language=Macaulay2output]
`\underline{\tt i6}` : gamma = makeGWClass matrix(RR, {{3,0,0},{0,-4,0},{0,0,7}});
`\underline{\tt i7}` : getSignature gamma
`\underline{\tt o7}` = `$1$`
\end{lstlisting}

Over finite fields, the \emph{discriminant}, which is the determinant of any Gram matrix representative (valued in square classes), and the rank jointly classify symmetric bilinear forms.

Over the rational numbers, the classification of symmetric bilinear forms is more complicated. The isomorphism class of a non-degenerate form $\beta$ is determined by its rank, discriminant $d_\beta\in \mathbb{Q^*}/(\mathbb{Q^*})^2$, signature (as a form over $\mathbb{R}$), and Hasse--Witt invariants at all finite primes~$p$. For a prime~$p$, the Hasse--Witt invariant \cite[III5, p.79]{MH} is defined as follows.
\begin{definition}\label{def:labelname} 
Given any form $\beta \cong \langle a_1, \dots, a_n \rangle \in \GW(\mathbb{Q})$, its \emph{Hasse--Witt invariant} $\varepsilon_p(\beta)$ at a prime~$p$
is the product
\[
\prod_{i<j} (a_i,a_j)_p,
\]
where $( -,- )_p$ denotes the \emph{Hilbert symbol}
\[
(a,b)_p \coloneqq \begin{cases} 
1 & \text{if } z^2 = ax^2 + by^2 \text{ has a nonzero solution in } \mathbb{Q}_p, \\
-1 & \text{otherwise.} \end{cases}
\]
\end{definition}
We can compute the Hilbert symbol $(a,b)_p$ via \texttt{getHilbertSymbol(a,b,p)} and the Hasse--Witt invariant of a form \texttt{beta} at~$p$ by \texttt{getHasseWittInvariant(beta, p)}. If $p$ is an odd prime and the $p$-adic valuations
$\nu_p(a)$ and $\nu_p(b)$ are even, then $(a,b)_p=1$. Thus, 
$\varepsilon_p(\beta)$ is~$1$ for almost all primes $p$ and only needs to be computed for $p=2$ and odd primes~$p$ with $\nu_p(d_\beta)$ odd.

These methods together form one of our core methods \texttt{isIsomorphicForm}, which is a Boolean-valued method that determines whether two symmetric bilinear forms are isomorphic. This is done by reference to the relevant invariants over $\mathbb{C}$, $\mathbb{R}$, $\mathbb{F}_q$ (for $q$ odd), or~$\mathbb{Q}$. 

\subsection{Decomposing forms}\label{subsec: decomposing forms}
\emph{Witt's Decomposition Theorem} (see \cite[I.4.1]{Lam}) implies that any non-degenerate symmetric bilinear form decomposes into an anisotropic part and an isotropic part that is a sum of hyperbolic forms. This decomposition is crucial in simplifying an element of $\GW(k)$. While this decomposition is fairly routine over $\mathbb{C}$, $\mathbb{R}$, and~$\mathbb{F}_q$, to decompose forms over~$\mathbb{Q}$ we must implement existing algorithms from the literature. An important mathematical stepping stone is the following local-to-global principle for isotropy, a reference for which is \cite[VI.3.1]{Lam}.

\begin{theorem}[Hasse--Minkowski Principle]\label{thm:lgp-isotropy}
A form $\beta\in \GW(\mathbb{Q})$ is isotropic if and only if it is isotropic over~$\mathbb{R}$ and over~$\mathbb{Q}_p$ for all primes~$p$.
\end{theorem}

Our method \texttt{getAnisotropicDimensionQQp}, an implementation of \cite[Algorithm~8]{KC}, determines the dimension of the anisotropic part of a symmetric bilinear form over~$\mathbb{Q}_p$. The method \texttt{getAnisotropicDimension} returns the anisotropic dimension of a form defined over the real numbers, the complex numbers, a finite field, or the rational numbers.

Given a form, we can therefore decompose it as
\[
\beta \cong \beta_a \oplus n \mathbb{H},
\]
where $\beta_a$ is anisotropic, $\mathbb{H}$ denotes the hyperbolic form $\langle 1,-1 \rangle$, and $n$ is the \emph{Witt index} (implemented as \texttt{getWittIndex}).

The Boolean-valued method \texttt{isAnisotropic} returns whether a form is anisotropic; the method \texttt{isIsotropic} is its negation.

\begin{lstlisting}[language=Macaulay2output]
`\underline{\tt i8}` : alpha = makeDiagonalForm(QQ, (1,2,-3));
`\underline{\tt i9}` : isAnisotropic alpha
`\underline{\tt o9}` = `$\texttt{false}$`
`\underline{\tt i10}` : isIsotropic alpha
`\underline{\tt o10}` = `$\texttt{true}$`
\end{lstlisting}

Over~$\mathbb{Q}$, the computation of the anisotropic part of~$\beta$ is carried out inductively by reduction of the anisotropic dimension of~$\beta$, following recently published algorithms for quadratic forms over number fields by Koprowski and Rothkegel~\cite{KR}. The anisotropic part of a form can be computed via \texttt{getAnisotropicPart}.

\begin{lstlisting}[language=Macaulay2output]
`\underline{\tt i11}` : beta = makeDiagonalForm(QQ, (3,-3,2,5,1,-9));
`\underline{\tt i12}` : getAnisotropicPart beta
`\underline{\tt o12}` = `$\left(\!\begin{array}{cc}
2&0\\
0&5
\end{array}\!\right)$`
`\underline{\tt o12}` : `$\texttt{GrothendieckWittClass}$`
\end{lstlisting}

A quick string reading off the decomposition of a form can be obtained by running the \texttt{getSumDecompositionString} method.
\begin{lstlisting}[language=Macaulay2output]
`\underline{\tt i13}` : getSumDecompositionString beta
`\underline{\tt o13}` = `\ttfamily 2H + {\char 60}2{\char 62} + {\char 60}5{\char 62}`
\end{lstlisting}

\section{\texorpdfstring{$\A^1$}{A\textonesuperior}-Brouwer degrees}\label{sec:A1Degs}

For the symbolic computations in this section, let $k$ be an exact field\footnote{An \emph{exact field} is a field whose elements are represented exactly by \emph{Macaulay2}, e.g., $\mathbb{Q}$ or $\mathbb{F}_q$.} of characteristic not equal to~$2$. The methods for computing $\mathbb{A}^{1}$-Brouwer degrees only work for polynomials with isolated zeros~\cite[Theorem~1.2]{BMP}.

In~\cite{BMP}, the authors show that the local and global $\mathbb{A}^{1}$-Brouwer degrees of an endomorphism of affine space with isolated zeros can be expressed in terms of a bilinear form associated to the B\'{e}zoutian of the endomorphism.

More explicitly, for $f_{i}\in k[x_{1},\dots,x_{n}]$, suppose $f=(f_{1},\dots,f_{n})\colon\mathbb{A}^{n}_{k}\to\mathbb{A}^{n}_{k}$ has isolated zeros. Introducing new variables $(X_{1},\dots,X_{n})$ and $(Y_{1},\dots, Y_{n})$, we can construct the matrix~$\Delta$ with entries 
\[
\Delta_{i,j}=\frac{f_{i}(Y_{1},\dots,Y_{j-1},X_{j},\dots,X_{n})-f_{i}(Y_{1},\dots,Y_{j},X_{j+1},\dots,X_{n})}{X_{j}-Y_{j}}.
\]
One can think of this matrix~$\Delta$ as a Jacobian of formal derivatives. Define $Q(f)=k[x_{1},\dots,x_{n}]/(f_{1},\dots,f_{n})$ and $Q_{p}(f)=k[x_{1},\dots,x_{n}]_{\mathfrak{m}}/(f_{1},\dots,f_{n})$ for $\mathfrak{m}$~the maximal ideal of a closed point~$p$ in the preimage of~$0$. The \emph{B\'{e}zoutian} of~$f$ is defined to be the image of $\det(\Delta)$ in the algebra $Q(f)\otimes Q(f)$ (respectively, in the local algebra $Q_{p}(f)\otimes Q_{p}(f)$). Given $a_{1},\dots,a_{m}$ a $k$-linear basis of $Q(f)$ (resp., $Q_{p}(f)$), there are $b_{i,j}\in k$ such that 
\[
\det(\Delta)=\sum_{1\leq i\leq j\leq m}b_{i,j}(a_{i}\otimes a_{j})
\]
in $Q(f)\otimes Q(f)$ (resp., $Q_{p}(f)\otimes Q_{p}(f)$). The \emph{B\'{e}zoutian bilinear form}, the symmetric bilinear form with Gram matrix given by the $b_{i,j}$, gives the global (resp., local) $\mathbb{A}^{1}$-degree \cite[Theorem~1.2]{BMP}. 

In the case of the global $\mathbb{A}^{1}$-degree, a theorem of Macaulay tells us that a $k$-basis of the algebra $Q(f)$ is given by the standard monomials (see \cite[Proposition~2.1]{BS}). In the case of the local $\mathbb{A}^{1}$-degree, a $k$-basis for the local ring can be calculated via the quotient of $k[x_{1},\dots,x_{n}]$ by the saturation of $I=(f_{1},\dots,f_{n})$ at~$\mathfrak{m}$. 
\begin{proposition}[{\cite[Proposition~2.5]{BS}}]\label{prop: local basis}
The natural map $x_i\mapsto x_i$ defines an isomorphism of rings
\begin{equation}\label{isomorph}
k[x_{1},\dots,x_{n}]_{\mathfrak{m}}/I\cong k[x_{1},\dots,x_{n}]/(I:(I:\mathfrak{m}^{\infty})),
\end{equation}
where $I$ is an ideal of $k[x_{1},\dots,x_{n}]$ and $(I:(I:\mathfrak{m}^{\infty}))$ is the quotient of~$I$ by the saturation of~$I$ at~$\mathfrak{m}$. 
\end{proposition}
The 
\texttt{getLocalAlgebraBasis(I, m)} method uses this isomorphism to find a basis of $Q_p(f)$. It 
determines
a $k$-basis of the right side of Equation~\eqref{isomorph} as a $k$-vector space. Proposition~\ref{prop: local basis} then gives a $k$-basis of $Q_p(f)$.

These methods for computing $k$-bases for $Q(f)$ and $Q_{p}(f)$ allow us to algorithmically implement techniques to compute the global and local $\mathbb{A}^{1}$-degrees (see also \cite[Section~5A]{BMP}). 

\subsection{A univariate polynomial}\label{subsec: univariate} 
A univariate polynomial over a field~$k$ defines an endomorphism of affine space $\mathbb{A}^{1}_{k}\to\mathbb{A}^{1}_{k}$. Consider the endomorphism $f\colon\mathbb{A}^{1}_{\mathbb{Q}}\to\mathbb{A}^{1}_{\mathbb{Q}}$ defined by 
\[
f(x)=(x^{2}+x+1)(x-3)(x+2).
\]
We can compute the global degree.

\begin{lstlisting}[language=Macaulay2output]
`\underline{\tt i14}` : R = QQ[x];
`\underline{\tt i15}` : f = {x^4 - 6*x^2 - 7*x - 6};
`\underline{\tt i16}` : alpha = getGlobalA1Degree f
`\underline{\tt o16}` = `$\left(\!\begin{array}{cccc}
-7&-6&0&1\\
-6&0&1&0\\
0&1&0&0\\
1&0&0&0
\end{array}\!\right)$`
`\underline{\tt o16}` : `$\texttt{GrothendieckWittClass}$`
\end{lstlisting}
We can also compute the local degrees at the ideals $(x^{2}+x+1)$, $(x-3)$, and $(x+2)$, respectively.

\begin{lstlisting}[language=Macaulay2output]
`\underline{\tt i17}` : I1 = ideal(x^2 + x + 1);
`\underline{\tt o17}` : `$\texttt{Ideal}$ of $R$`
`\underline{\tt i18}` : alpha1 = getLocalA1Degree(f, I1)
`\underline{\tt o18}` = `$\left(\!\begin{array}{cc}
-5&-7\\
-7&-2
\end{array}\!\right)$`
`\underline{\tt o18}` : `$\texttt{GrothendieckWittClass}$`
`\underline{\tt i19}` : I2 = ideal(x - 3)
`\underline{\tt o19}` = `$\texttt{ideal}\left(x-3\right)$`
`\underline{\tt o19}` : `$\texttt{Ideal}$ of $R$`
`\underline{\tt i20}` : alpha2 = getLocalA1Degree(f, I2)
`\underline{\tt o20}` = `$\left(\!\begin{array}{c}
65
\end{array}\!\right)$`
`\underline{\tt o20}` : `$\texttt{GrothendieckWittClass}$`
`\underline{\tt i21}` : I3 = ideal(x + 2);
`\underline{\tt o21}` : `$\texttt{Ideal}$ of $R$`
`\underline{\tt i22}` : alpha3 = getLocalA1Degree(f, I3)
`\underline{\tt o22}` = `$\left(\!\begin{array}{c}
-15
\end{array}\!\right)$`
`\underline{\tt o22}` : `$\texttt{GrothendieckWittClass}$`
\end{lstlisting}
We can then use the \texttt{isIsomorphicForm} method (see also Section~\ref{subsec: isomorphisms of forms}) to verify that the local $\mathbb{A}^{1}$-degrees sum to the global $\mathbb{A}^{1}$-degree.

\begin{lstlisting}[language=Macaulay2output]
`\underline{\tt i23}` : alpha' = addGW(alpha1,addGW(alpha2,alpha3));
`\underline{\tt i24}` : isIsomorphicForm(alpha,alpha')
`\underline{\tt o24}` = `$\texttt{true}$`
\end{lstlisting}

Consider the graph of $f(x)$. 
\begin{center}
    \begin{tikzpicture}
    \begin{axis} [axis lines=center]
    \addplot [domain=-3:3.7, smooth, thick] { x^4 - 6*x^2 - 7*x - 6 };
  \end{axis}
\end{tikzpicture}
\end{center}
Following~\cite[Theorem 1.2]{McKean}, $\mathbb{A}^{1}$-degrees can be understood as enriched intersection numbers, determined by the signed volume of the parallelepiped spanned by the gradient vectors of the hypersurfaces at the intersection point. In the one-dimensional case, considering the normal vectors, we can interpret $\alpha_{2}=\langle65\rangle$, the local $\mathbb{A}^{1}$-degree at $(x-3)$, and $\alpha_{3}=\langle-15\rangle$, the local $\mathbb{A}^{1}$-degree at $(x+2)$, as signs of the derivative at these points. 

\subsection{The Euler characteristic of the Grassmannian of lines in \texorpdfstring{$\mathbb{P}^{3}$}{P\textthreesuperior}}\label{subsec: Grassmannian} 
For $k$ a field of characteristic not~$2$, let $\Gr_{k}(2,4)$ be the Grassmannian of lines in~$\mathbb{P}^{3}_{k}$. Following \cite[Example~8.2]{BMP}, we can compute the $\mathbb{A}^{1}$-Euler characteristic of the Grassmannian over $k=\mathbb{F}_{27}$ as the $\mathbb{A}^{1}$-degree of the section $\sigma\colon\mathbb{A}^{4}_{\mathbb{F}_{27}}\to\mathbb{A}^{4}_{\mathbb{F}_{27}}$ defined by\footnote{%
There is a small error in the definition of $\sigma$ in \cite[Example~8.2]{BMP}. The second and third component functions of $\sigma$ should be swapped in order to agree with the ordered basis induced on the tangent bundle of $\Gr_k(2,4)$ as in \cite[Proposition~45]{CubicKW}. By \cite[Example~6.3]{BMP}, the overall computation is only affected by a sign.
}
\[
(x_1,x_2,x_3,x_4)\mapsto(
x_{2}-x_{1}x_{3},
1-x_{1}x_{4},
x_{4}-x_{1}-x_{3}^{2},
-x_{2}-x_{3}x_{4}).
\]
We compute the $\mathbb{A}^{1}$-Euler characteristic as follows.

\begin{lstlisting}[language=Macaulay2output]
`\underline{\tt i25}` : k = GF(27);
`\underline{\tt i26}` : x = symbol x;
`\underline{\tt i27}` : R = k[x_1,x_2,x_3,x_4];
`\underline{\tt i28}` : f = {x_2 - x_1*x_3, 1 - x_1*x_4, x_4 - x_1 - x_3^2, -x_2 - x_3*x_4};
`\underline{\tt i29}` : beta = getGlobalA1Degree f
`\underline{\tt o29}` = `$\left(\!\begin{array}{cccccc}
0&0&0&0&0&1\\
0&1&0&0&0&0\\
0&0&0&-1&0&0\\
0&0&-1&0&0&0\\
0&0&0&0&1&0\\
1&0&0&0&0&0
\end{array}\!\right)$`
`\underline{\tt o29}` : `$\texttt{GrothendieckWittClass}$`
\end{lstlisting}
We can subsequently use the \texttt{getSumDecompositionString} method to decompose the symmetric bilinear form~$\beta$.

\begin{lstlisting}[language=Macaulay2output]
`\underline{\tt i30}` : getSumDecompositionString beta
`\underline{\tt o30}` = `\ttfamily 2H + {\char 60}1{\char 62} + {\char 60}1{\char 62}`
\end{lstlisting}
Our computation agrees with the result given in \cite[Example~8.2]{BMP} and shows \[\chi(\Gr_{\mathbb{F}_{27}}(2,4))=2\mathbb{H}+\langle 1\rangle + \langle 1\rangle.\]

\subsection{Local geometry of some lines on the Fermat cubic surface} 
In their pioneering paper~\cite{CubicKW}, Kass and Wickelgren give a Grothendieck--Witt class-valued count of the number of lines on a smooth cubic surface, providing an interpretation of the local $\mathbb{A}^{1}$-degree as the topological type of the line. To illustrate some features of the $\mathtt{A1BrouwerDegrees}$ package, we use it to compute the topological type of some lines on the Fermat cubic surface.

Let $k$ be a field, and let $\{e_{1},e_{2},e_{3},e_{4}\}$ be the standard basis for~$k^{4}$. By \cite[Lemma~45]{CubicKW}, we can define local coordinates on $\Spec(k[y_{1},y_{2},y_{3},y_{4}])\cong\mathbb{A}_k^{4}$ around the point of $\Gr_k(2,4)$ defined by the span of $\{e_{3},e_{4}\}$ such that $y_{1},y_{2},y_{3},y_{4}$ corresponds to the span of $\{\widetilde{e}_3,\widetilde{e}_{4}\}$, where 
\[
\widetilde{e}_i = \begin{cases}
    e_{i} & \text{for } i\in\{1,2\}, \\
    e_{1}y_{1} + e_{2}y_{2} + e_{3} & \text{for } i=3, \\
    e_{1}y_{3} + e_{2}y_{4} + e_{4} & \text{for } i=4.
\end{cases}
\]
Letting $\mathcal{S}$ denote the tautological bundle over $\Gr_k(2,4)$, the above coordinates provide a trivialization of the vector bundle $\Sym^{3}\mathcal{S}^{\vee}$ over the open affine subvariety
\[
U\cong\Spec(k[y_{1},y_{2},y_{3},y_{4}])\subseteq\Gr_k(2,4).
\]
A cubic surface~$X$ defines a section $\sigma_{X}\rvert_{U}\colon U\to\Sym^{3}\mathcal{S}^{\vee}\rvert_{U}$ that vanishes on the lines on~$X$ that, when treated as affine two-dimensional subspaces of~$k^{4}$, contain $e_{3}$ and~$e_{4}$ in their span. 

Let us consider the Fermat cubic surface defined by the homogeneous cubic equation $x_{0}^{3}+x_{1}^{3}+x_{2}^{3}+x_{3}^{3}$. That is, 
\[
X=\left\{[x_{0}:x_{1}:x_{2}:x_{3}]\in\mathbb{P}_k^{3}:x_{0}^{3}+x_{1}^{3}+x_{2}^{3}+x_{3}^{3}=0\right\}\subseteq\mathbb{P}_k^{3}.
\]
Working over~$\mathbb{Q}$, the lines on~$X$ are all defined over the cyclotomic extension $\mathbb{Q}(\zeta)$ for~$\zeta$ a primitive third root of unity. We can explicitly compute the $27$~lines as
\[
[s:t:-\zeta^{i}t:-\zeta^{j}s], [s:t:-\zeta^{i}s:-\zeta^{j}t], [s:-\zeta^{i}s: t:-\zeta^{j}t]
\]
for $0\leq i,j\leq 2$ and $[s:t]\in\mathbb{P}_{\mathbb{Q}}^{1}$. Note that there are only $18$~lines containing $e_{3}$ and~$e_{4}$ in their span. We thus expect the section to vanish at $18$~points. Applying Pauli's computation of sections of $\sigma_{X}\rvert_{U}\colon U\to\Sym^{3}\mathcal{S}^{\vee}\rvert_{U}$ in \cite[\S 2.2, Remark 2.7]{Pauli}, our section is of the form $\sigma_{X}\rvert_{U}=(f_{1},f_{2},f_{3},f_{4})$ where
\begin{align*}
    f_{1}(y_{1},y_{2},y_{3},y_{4}) &= y_{1}^{3}+y_{3}^{3}+1 \\
    f_{2}(y_{1},y_{2},y_{3},y_{4}) &= 3y_{1}^{2}y_{2}+3y_{3}^{2}y_{4}\\
    f_{3}(y_{1},y_{2},y_{3},y_{4}) &= 3y_{1}y_{2}^{2}+3y_{3}y_{4}^{2}\\
    f_{4}(y_{1},y_{2},y_{3},y_{4}) &= y_{2}^{3}+y_{4}^{3}+1.
\end{align*}
We compute the global $\mathbb{A}^{1}$-degree, which is rank~$18$, as expected. 

\begin{lstlisting}[language=Macaulay2output]
`\underline{\tt i31}` : R = QQ[y_1,y_2,y_3,y_4];
`\underline{\tt i32}` : f = {y_1^3 + y_3^3 + 1, 
          3*y_1^2*y_2 + 3*y_3^2*y_4, 
          3*y_1*y_2^2 + 3*y_3*y_4^2, 
          y_2^3 + y_4^3 + 1};
`\underline{\tt i33}` : alpha = getGlobalA1Degree f;
`\underline{\tt i34}` : getSumDecompositionString alpha
`\underline{\tt o34}` = `\ttfamily 8H + {\char 60}1{\char 62} + {\char 60}1{\char 62}`
\end{lstlisting}
To compute the local degree, we find an isolated zero using the \texttt{minimalPrimes} method of \emph{Macaulay2}.

\begin{lstlisting}[language=Macaulay2output]
`\underline{\tt i35}` : I = (minimalPrimes ideal f)_0
`\underline{\tt o35}` = `$\texttt{ideal}{}\left(y_{4},\,y_{3}+1,\,y_{2}+1,\,y_{1}\right)$`
`\underline{\tt o35}` : `$\texttt{Ideal}$ of $R$`
\end{lstlisting}
We then compute the local $\mathbb{A}^{1}$-degree at this point.

\begin{lstlisting}[language=Macaulay2output]
`\underline{\tt i36}` : beta = getLocalA1Degree(f, I)
`\underline{\tt o36}` = `$\left(\!\begin{array}{c}
81
\end{array}\!\right)$`
`\underline{\tt o36}` : `$\texttt{GrothendieckWittClass}$`
`\underline{\tt i37}` : getSumDecomposition beta
`\underline{\tt o37}` = `$\left(\!\begin{array}{c}
1
\end{array}\!\right)$`
`\underline{\tt o37}` : `$\texttt{GrothendieckWittClass}$`
\end{lstlisting}
At the other nine minimal primes, the same calculation gives the local degrees as one copy of $\langle 1\rangle$, six copies of $\langle 3, -1\rangle$ and two copies of $\langle2, -6\rangle$. As 
\[
6\langle3, -1\rangle+2\langle2,-6\rangle\cong 8\mathbb{H},
\]
the global $\mathbb{A}^1$-degree equals the sum of the local $\mathbb{A}^1$-degrees.

The local computation above indicates that the line spanned by $\{-e_{2}+e_{3},-e_{1}+e_{4}\}$ on the Fermat cubic surface is a hyperbolic line. We briefly show this agrees with the type as defined in \cite[Definition~9]{CubicKW}. 

By \cite[Proposition~14]{CubicKW}, the local type of the line is equal to the resultant of the partial derivatives of the equation of the Fermat cubic surface restricted to the line. Letting $z_{1},z_{2},z_{3},z_{4}$ be the dual basis to $e_{1},e_{2},e_{3},e_{4}$ defined above, we can write the equation of the Fermat surface in terms of the dual basis via the change of basis $z_{1}\mapsto z_{1}+z_{4}, z_{2}\mapsto z_{2}+z_{3}, z_{3}\mapsto -z_{3}, z_{4}\mapsto -z_{4}$ so that the line is spanned by $e_{3}$ and $e_{4}$. 

\begin{lstlisting}[language=Macaulay2output]
`\underline{\tt i38}` : needsPackage "Resultants";
`\underline{\tt i39}` : R = QQ[z_1,z_2][z_3,z_4];
`\underline{\tt i40}` : fermat = (z_1 + z_4)^3 + (z_2 + z_3)^3 - z_3^3 - z_4^3;
\end{lstlisting}
We compute the restriction of the partial derivatives of the defining equation of the Fermat cubic surface to the surface with respect to the dual basis $z_{1},z_{2}$. 

\begin{lstlisting}[language=Macaulay2output]
`\underline{\tt i41}` : g1 = sub(diff(z_1, fermat), {z_1 => 0, z_2 => 0});
`\underline{\tt i42}` : g2 = sub(diff(z_2, fermat), {z_1 => 0, z_2 => 0});
\end{lstlisting}
We then compute the resultant of these polynomials and consider it as a quadratic form over~$\mathbb{Q}$ in order to agree with the computation of the local index over~$\mathbb{Q}$.
\begin{lstlisting}[language=Macaulay2output]
`\underline{\tt i43}` : line_type = makeDiagonalForm(QQ, lift(resultant {g1,g2}, QQ))
`\underline{\tt o43}` = `$\left(\!\begin{array}{c}
81
\end{array}\!\right)$`
`\underline{\tt o43}` : `$\texttt{GrothendieckWittClass}$`
`\underline{\tt i44}` : isIsomorphicForm(line_type,beta)
`\underline{\tt o44}` = `$\texttt{true}$`
\end{lstlisting}
Thus this computation agrees with the local $\A^1$-degree of the associated section of $\Sym^3 \mathcal{S}^\ast$ as computed above.

\section*{Acknowledgements}

We would like to thank the organizers of the 2023 \emph{Macaulay2} Workshop for the opportunity to work together developing this package. The second-named author would like to thank Sabrina Pauli for helpful conversations about this work. We would also like to thank the anonymous referees for detailed comments that helped to improve both this paper and the code. The second-named author is supported by an NSF Postdoctoral Research Fellowship (DMS-2303242), and the sixth-named author is partially supported by an NSF Standard Grant of Frank Sottile (DMS-2201005).

\bibliographystyle{abbrv}
\bibliography{jsag}

\begin{thebibliography}{10}

\bibitem{AK}
D.~Agostini and M.~Kummer.
\newblock Secant varieties of curves, {U}lrich bundles and the arithmetic
  writhe, 2023.
\newblock \href{https://arxiv.org/abs/2307.07543}{arXiv:2307.07543} [math.AG].

\bibitem{BW}
T.~Bachmann and K.~Wickelgren.
\newblock Euler classes: six-functors formalism, dualities, integrality and
  linear subspaces of complete intersections.
\newblock {\em J. Inst. Math. Jussieu}, 22(2):681--746, 2023.

\bibitem{Braz-expos}
T.~Brazelton.
\newblock An introduction to {$\mathbb{A}^1$}-enumerative geometry.
\newblock In {\em Homotopy theory and arithmetic geometry---motivic and
  {D}iophantine aspects}, volume 2292 of {\em Lecture Notes in Math.}, pages
  11--47. Springer, Cham, 2021.

\bibitem{BMP}
T.~Brazelton, S.~McKean, and S.~Pauli.
\newblock B\'{e}zoutians and the {$\mathbb A^1$}-degree.
\newblock {\em Algebra Number Theory}, 17(11):1985--2012, 2023.

\bibitem{EL}
D.~Eisenbud and H.~I. Levine.
\newblock An algebraic formula for the degree of a {$C^{\infty }$} map germ.
\newblock {\em Ann. of Math. (2)}, 106(1):19--44, 1977.
\newblock With an appendix by Bernard Teissier, ``Sur une in\'{e}galit\'{e} \`a
  la Minkowski pour les multiplicit\'{e}s''.

\bibitem{M2}
D.~R. Grayson and M.~E. Stillman.
\newblock Macaulay2, a software system for research in algebraic geometry.
\newblock Available at \url{http://www2.macaulay2.com}.

\bibitem{Khim}
G.~N. Him\v{s}ia\v{s}vili.
\newblock The local degree of a smooth mapping.
\newblock {\em Sakharth. SSR Mecn. Akad. Moambe}, 85(2):309--312, 1977.

\bibitem{JPP}
A.~{Jaramillo Puentes} and S.~Pauli.
\newblock A quadratically enriched correspondence theorem, 2024.
\newblock \href{https://arxiv.org/abs/2309.11706}{arXiv:2309.11706} [math.AG].

\bibitem{KSLW}
J.~L. Kass, M.~Levine, J.~P. Solomon, and K.~Wickelgren.
\newblock A quadratically enriched count of rational curves, 2023.
\newblock \href{https://arxiv.org/abs/2307.01936}{arXiv:2307.01936} [math.AG].

\bibitem{KW-EKL}
J.~L. Kass and K.~Wickelgren.
\newblock The class of {E}isenbud-{K}himshiashvili-{L}evine is the local
  {$\mathbf{A}^1$}-{B}rouwer degree.
\newblock {\em Duke Math. J.}, 168(3):429--469, 2019.

\bibitem{CubicKW}
J.~L. Kass and K.~Wickelgren.
\newblock An arithmetic count of the lines on a smooth cubic surface.
\newblock {\em Compos. Math.}, 157(4):677--709, 2021.

\bibitem{KC}
P.~Koprowski and A.~Czoga{\l}a.
\newblock Computing with quadratic forms over number fields.
\newblock {\em J. Symbolic Comput.}, 89:129--145, 2018.

\bibitem{KR}
P.~Koprowski and B.~Rothkegel.
\newblock The anisotropic part of a quadratic form over a number field.
\newblock {\em J. Symbolic Comput.}, 115:39--52, 2023.

\bibitem{Lam}
T.~Y. Lam.
\newblock {\em Introduction to quadratic forms over fields}, volume~67 of {\em
  Graduate Studies in Mathematics}.
\newblock American Mathematical Society, Providence, RI, 2005.

\bibitem{Levine}
M.~Levine.
\newblock Aspects of enumerative geometry with quadratic forms.
\newblock {\em Doc. Math.}, 25:2179--2239, 2020.

\bibitem{McKean}
S.~McKean.
\newblock An arithmetic enrichment of {B}\'{e}zout's {T}heorem.
\newblock {\em Math. Ann.}, 379(1--2):633--660, 2021.

\bibitem{MH}
J.~Milnor and D.~Husemoller.
\newblock {\em Symmetric bilinear forms}, volume Band 73 of {\em Ergebnisse der
  Mathematik und ihrer Grenzgebiete [Results in Mathematics and Related
  Areas]}.
\newblock Springer-Verlag, New York-Heidelberg, 1973.

\bibitem{Morel}
F.~Morel.
\newblock {\em {$\mathbb{A}^1$}-algebraic topology over a field}, volume 2052
  of {\em Lecture Notes in Mathematics}.
\newblock Springer, Heidelberg, 2012.

\bibitem{Pauli}
S.~Pauli.
\newblock Computing {$\mathbb{A}^1$}-{E}uler numbers with {M}acaulay2.
\newblock {\em Res. Math. Sci.}, 10(3):Paper No. 26, 14 pp., 2023.

\bibitem{Pauli-expos}
S.~Pauli and K.~Wickelgren.
\newblock Applications to {$\mathbb{A}^1$}-enumerative geometry of the
  {$\mathbb{A}^1$}-degree.
\newblock {\em Res. Math. Sci.}, 8(2):Paper No. 24, 29 pp., 2021.

\bibitem{BS}
B.~Sturmfels.
\newblock {\em Solving systems of polynomial equations}, volume~97 of {\em CBMS
  Regional Conference Series in Mathematics}.
\newblock Published for the Conference Board of the Mathematical Sciences,
  Washington, DC; by the American Mathematical Society, Providence, RI, 2002.

\end{thebibliography}

\end{document}